\documentclass{amsart}

\usepackage{amssymb, mathrsfs}
\usepackage[all]{xy}
\usepackage{latexsym}

\usepackage{amsmath}
\usepackage{amsthm}

\theoremstyle{plain}
\newtheorem{theorem}{Theorem}
\newtheorem{definition}[theorem]{Definition} 
\newtheorem{lemma}[theorem]{Lemma}
\newtheorem{proposition}[theorem]{Proposition}
\newtheorem{corollary}[theorem]{Corollary}
\newtheorem{remark}[theorem]{Remark}

\theoremstyle{definition}
\newtheorem{example}[theorem]{Example}

\providecommand{\Hom}{\mathop{\rm Hom}\nolimits}
\providecommand{\dim} {\mathop{\rm dim}}

\providecommand{\rad} {\mathop{\rm rad}\nolimits}
\providecommand{\HH} {\rm {HH}}
\providecommand{\SH}{\rm {SH}}
\providecommand{\Ext} {\mathop{\rm Ext}\nolimits}
\providecommand{\Ker} {\mathop{\rm Ker}\nolimits}

\providecommand{\Ia}{{\mathcal I}}

\begin{document}

\title[Hochschild cohomology via incidence algebras]{Hochschild cohomology via incidence algebras}

\author{Mar\'\i a Julia Redondo}
\address{Maria Julia Redondo \\ 
Instituto de Matem\'atica, Universidad Nacional del Sur \\ 
Av. Alem 1253, (8000) Bah\'\i a Blanca, Argentina.} 



\begin{abstract}
Given an algebra $A$ we associate an incidence algebra $A(\Sigma)$ and compare their Hochschild cohomology groups.
\end{abstract}

\maketitle

\section {Introduction}

The purpose of this paper is to study the connection between the Hochschild cohomology groups of $A$ and of an incidence algebra associated to each presentation $(Q,I)$ of $A$.

Let $A$ be an associative, finite dimensional algebra over an algebraically closed field $k$.  It is well known that if $A$ is basic and connected, then there exists a unique finite quiver $Q$ such that $A \cong kQ/I$, where $kQ$ is the path algebra of $Q$ and $I$ is an admissible two-sided ideal of $kQ$.  The pair $(Q,I)$  is called a \textit{presentation} of $A$.

For each presentation $(Q,I)$ of $A$, one can define its \textit{fundamental group} $\pi_1(Q,I)$.  This group is not an invariant of the algebra, since different presentations may yield the same algebra, with different fundamental groups associated.  In~\cite{AP,PS} it is shown that this group is related to the first Hochschild cohomology group of $A$: for any presentation $(Q,I)$ of $A$ there is an injective morphism of abelian groups $\Hom(\pi_1(Q,I),k^{+}) \to {\HH}^1(A)$, where $k^{+}$ denotes the underlying additive group of the field $k$.

If $A$ is an \textit{incidence algebra}, the fundamental group is an invariant of the algebra, which will be denoted by $\pi_1(A)$, and  $\Hom(\pi_1(A),k^{+}) \simeq {\HH}^1(A)$.  Moreover, there is a simplicial complex $\vert A \vert$  such that the fundamental group $\pi_1(A)$ and the (co)homology groups of $A$ are respectively isomorphic to the fundamental group and the (co)homology groups of $\vert A \vert$, see~\cite{B,Re1}.

In~\cite{Re2} it is shown that for each presentation $(Q,I)$ of an algebra $A$ there is an associated incidence algebra $A(\Sigma)$, and the corresponding fundamental groups are related by the following short exact sequence of groups
$$ 1 \to H \to \pi_1(Q,I) \to \pi_1(A(\Sigma)) \to 1.$$
Moreover, $H$ is explicitly described in~\cite{Re2} by generators and relations.

Even though the computation of the Hochschild cohomology groups ${\HH}^i(A)$ is rather complicated, some approaches have been successful when the algebra $A$ is
given by a quiver with relations. For instance, explicit formula for the dimensions of ${\HH}^i(A)$ in terms of
those combinatorial data have been found in~\cite{blm,c1,C2,crs,cs,H}. 

In~\cite{IZ} Igusa and Zacharia give a combinatorial algorithm to find an upper bound for the cohomological dimension of ${\HH}^i({\Ia}(\Sigma))$, where ${\Ia}(\Sigma)$ is the incidence algebra associated to the poset $\Sigma$. They show how to construct the so-called {\it reduced subposet} $\overline{\Sigma}$ of $\Sigma$ which has the property that all elements $x \in \overline{\Sigma}$, neither minimal nor maximal elements, are such that $\{y \in \overline{\Sigma}: y > x\}$ has at least two minimal elements and $\{z \in \overline{\Sigma}: x > z \}$ has at least two maximal elements. The Hochschild cohomology groups are invariant under this construction.  Hence it is enough to compute them for incidence algebras associated to reduced posets. The Hochschild cohomology groups of some particular families of incidence algebras have been computed in~\cite{GR,GR1} and it is known that they vanish if the associated poset does not contain crowns, see~\cite{D,IZ}.

In this paper we are interested in comparing the Hochschild cohomology groups of an algebra $A$ with those of a particular incidence algebra $A(\Sigma)$ associated to a presentation $(Q,I)$ of $A$.  If the chosen presentation is homotopy coherent, see Definition~\ref{P1}, we define a morphism between the complexes computing these cohomology groups, which induces morphisms ${\HH}(\Phi^n): {\HH}^n(A(\Sigma_\nu)) \to {\HH}^n(A)$. Finally we find conditions for these morphisms to be injective, see Theorem~\ref{teorema}.

The paper is organized as follows.  In Section 2 we introduce all the necessary terminology and known results.  In Section 3 we recall the construction of the incidence algebra associated to any presentation of an algebra $A$.  In Section 4 we introduce two classes of presentations: the homotopy coherent and the right (left) compatible presentations, that will be essentially needed to prove the main results, presented in Section 5, concerning the relationship between the Hochschild cohomology groups of $A$ and those of an associated incidence algebra.  Finally in Section 6 we present several examples.

\section {Preliminaries}

\subsection{Quivers and relations}

Let $Q$ be a  finite quiver with a set of vertices $Q_0$, a set of arrows $Q_1$ and $s, t : Q_1 \to Q_0$ be the maps associating to each arrow $\alpha$ its source  $s(\alpha)$ and its target $t(\alpha)$.  A path $w$ of length $l$ is a sequence of $l$ arrows $\alpha_1 \dots \alpha_l$ such that $t(\alpha_i)=s(\alpha_{i+1})$.  We put $s(w)=s(\alpha_1)$ and $t(w)=t(\alpha_l)$. For any vertex $x$ we consider $e_x$ the trivial path of length zero and we put $s(e_x)=t(e_x)=x$.  A cycle is a non-trivial path $w$ such that $s(w)=t(w)$. 

The corresponding path algebra $kQ$ is the $k$-vector space with basis the set of paths in $Q$; the product on the basis elements is given by the concatenation of the sequences of arrows of the paths $w$ and $w'$ if they form a path (namely, if $t(w)=s(w')$) and zero otherwise.  Vertices form a complete set of orthogonal idempotents. Let $F$ be the two-sided ideal of $kQ$ generated by the arrows of $Q$. A two-sided ideal $I$ is said to be \textit{admissible} if there exists an integer $m \geq 2$ such that $F^m \subseteq I \subseteq F^2$.  The pair $(Q,I)$  is called a \textit{bound quiver}.

It is well known that if $A$ is a basic, connected, finite dimensional algebra over an algebraically closed field $k$, then there exists a unique finite quiver $Q$ and a surjective morphism of $k$-algebras $\nu: kQ \to A$, which is not unique in general, with $I_\nu=\Ker \nu$ admissible.  The pair $(Q,I_\nu)$  is called a \textit{presentation} of $A$. We denote by $kQ(x,y)$ the subspace of $kQ$ with basis the set of paths from $x$ to $y$, $I_\nu(x,y)=I_\nu \cap kQ(x,y)$ and $A(x,y)=kQ(x,y)/I_\nu(x,y)$.

\subsection{Incidence algebras}\label{incidencia}

An incidence algebra ${\Ia}(\Sigma)$ is a subalgebra of the algebra $M_n(k)$ of square matrices over $k$ with elements $(x_{ij}) \in M_n(k)$ satisfying $x_{ij}=0$ if $i \not \geq j$, for some partial order $\ge$ defined in the {\it poset} (partially ordered set) $\Sigma= \{ 1, \dots , n\}$. 

Incidence algebras can equivalently be viewed as path algebras of quivers with relations in the following way. Let $Q$ be a finite quiver without oriented cycles and such that for each arrow $ x \stackrel{\alpha}{\rightarrow} y \in Q_1$ there is no oriented path other than $\alpha$ joining $x$ to $y$. These quivers are called {\it ordered}. 
The set $Q_0$ of vertices of $Q$ is then a finite poset $(\Sigma, \geq)$ as follows: $x \geq y$ if and only if there exists an oriented path from $x$ to $y$. Conversely, if $\Sigma$ is a finite poset, we construct a quiver $Q$ with the set of vertices $\Sigma$, and with an arrow from $x$ to $y$ if and only if $x>y$ and there is no $u \in Q_0$ such that $x>u>y$. In other words $Q$ is the Hasse diagram of the poset $\Sigma$. Clearly we obtain in this way an ordered quiver and a bijection between finite posets and ordered quivers.

Let us consider $kQ$ the path algebra of $Q$ and $I$ the {\it parallel ideal} of $kQ$, that is, $I$ is the two-sided ideal of $kQ$ generated by all the differences $\gamma - \delta$ where $\gamma$ and $\delta$ are 
{\it parallel paths} (that is, $\gamma$ and $\delta$ have the same starting and ending points). The algebra ${\Ia}(\Sigma) = kQ/I$ is the {\it incidence algebra} of the poset $\Sigma$ associated to the ordered quiver $Q$. 

\subsection{Fundamental group}

Let $(Q,I)$ be a connected bound quiver. For $x,y \in Q_0$, a relation $\rho = \sum_{i=1}^m \lambda_i w_i$ in $I(x,y)$ is called {\it minimal} if, for every non-empty proper subset $J \subset \{ 1,2, \dots,m\}$, we have $\sum_{j \in J} \lambda_j w_j \notin I(x,y)$. A relation $\rho$ is called \textit{monomial} if $m=1$, and \textit{binomial} if $m=2$.

For an arrow $\alpha \in Q_1$, we denote by $\alpha^{-1}$ its formal inverse and put $s(\alpha^{-1})=t(\alpha)$ and $t(\alpha^{-1})=s(\alpha)$. A walk from $x$ to $y$ in $Q$ is a formal composition  $\alpha_1^{\epsilon_1}\alpha_2^{\epsilon_2} \dots \alpha_t^{\epsilon_t}$  
(where $\alpha_i \in Q_1$, $\epsilon_i=^+_-1$ for $1\leq i\leq t$, and $t(\alpha_i^{\epsilon_i})=s(\alpha_{i+1}^{\epsilon_{i+1}})$) starting at $x$ and ending at $y$. 

Let $\approx$ be the smallest equivalence relation on the set of all walks in $Q$ such that:
\begin{itemize}
\item[a)] If $\alpha: x \to y$ is an arrow then $\alpha^{-1}\alpha \approx e_y$ and $\alpha \alpha^{-1} \approx e_x$;
\item[b)] If $\rho = \sum_{i=1}^m \lambda_iw_i$ is a minimal relation then $w_i \approx w_j$ for all $1 \leq i,j \leq m$;
\item[c)] If $u \approx v$ then $wuw' \approx wvw'$ whenever these compositions make sense.
\end{itemize}
Let $x_0 \in Q_0$ be arbitrary. The set $\pi_1(Q,I,x_0)$ of equivalence classes of all the closed walks starting and ending at $x_0$ has a group structure. Clearly the group $\pi_1(Q,I,x_0)$ does not depend on the choice of the base point $x_0$. We denote it simply by $\pi_1 (Q,I)$ and call it the {\it fundamental group} of $(Q,I)$, see~\cite{AP}.

\subsection{Hochschild cohomology: a convenient resolution for $A=kQ/I$}

We recall that the Hochschild cohomology groups ${\HH}^i(A)$ of an algebra $A$ are the groups $\Ext_{A^e}^i(A,A)$.
We refer the reader to~\cite{CE,H,R} for more general results.

To compute the Hochschild cohomology groups of $A=kQ/I$ we will use a convenient projective resolution of $A$ as $A$-bimodule given in~\cite {C}.
Let $E$ be the subalgebra of $A$ generated by the set of trivial paths $\{e_x \vert x \in Q_0 \}$; note that $E$ is semisimple and $A=E \oplus \rad A$ as $E$-bimodule.  Let $\rad A ^{\otimes n}$ denote the $n$-fold tensor product of $\rad A$ with itself over $E$, with $\rad A^{\otimes 0}=E$. The complex 
$$ \cdots  \to A \otimes_E \rad A^{\otimes 2} \otimes_E A \stackrel{b_2}{\to} A \otimes_E \rad A \otimes_E A \stackrel{b_1}{\to}  A \otimes_E A \stackrel{b_0}{\to} A \to 0$$
with
$$b_n(a_0 \otimes \dots \otimes a_{n+1}) =  \sum_{i=0}^{n} {(-1)}^i  a_0\otimes \dots \otimes a_i a_{i+1}
\otimes \dots \otimes a_{n+1}$$
is a projective resolution of $A$ as $A$-bimodule.
There is a natural isomorphism $\Hom_{A^e}(A \otimes_E \rad A^{\otimes n} \otimes_E A, A) \simeq \Hom_{E^e}(\rad A^{\otimes n}, A)$, 
and the corresponding boundary map $b^n: \Hom_{E^e}(\rad A^{\otimes n}, A) \to \Hom_{E^e}(\rad A^{\otimes n+1}, A)$ is given by
\begin{equation*}
\begin{split}
(b^0f)(a_0) & =a_0f(1)-f(1)a_0,\\
(b^nf)(a_0 \otimes \dots \otimes a_{n}) &= a_0f(a_1 \otimes \dots \otimes a_{n}) \\
& \quad + \sum_{i=1}^{n} {(-1)}^{i}f(a_0\otimes \dots \otimes a_{i-1} a_{i}
\otimes \dots \otimes a_{n}) \\
& \quad +{(-1)}^{n+1} f(a_0 \otimes \dots \otimes a_{n-1})a_{n}.
\end{split}
\end{equation*}

\subsection{Hochschild cohomology and simplicial cohomology} \label{posetincidencia}

Let $A={\Ia}(\Sigma)$ be an incidence algebra associated to a poset $\Sigma$.  The simplicial complex associated to $\Sigma$ is defined as follows: $SC_n=SC_n(\Sigma)$ is the $k$-vector space with basis the set $\{ s_0 > s_1 > \dots > s_n \vert s_i \in \Sigma\}$. The complex computing the cohomology groups $\SH^n(\Sigma,k)$ is the following
$$ 0 \to \Hom_k(SC_0,k) \stackrel{B^0}{\to} \Hom_k(SC_1,k)  \stackrel{B^1}{\to} \Hom_k(SC_2,k) \to \cdots$$
with 
$$(B^{n-1} f)(s_0> \dots > s_n) = \sum_{i=0}^n (-1)^i f(s_0> s_1 > \dots > \hat{s_i} > \dots > s_n).$$
In~\cite{GS,C} it was shown that $\SH^n(\Sigma,k)$ and ${\HH}^n(A)$ are isomorphic. Moreover, there is an explicit isomorphism between the complexes computing these cohomology groups, which is defined as follows: note that $A={\Ia}(\Sigma)$ is an incidence algebra, and we are considering tensor products over $E$, thus $$\rad A^{\otimes n}=  \oplus_{s_0> s_1 > \dots > s_n} A(s_0, s_1) \otimes_k A(s_1, s_2) \otimes_k \dots \otimes_k A(s_{n-1}, s_n),$$
with $\dim_k A(s_{i-1}, s_i)=1$ for all $i$ with $1\leq i \leq n$. Taking a basis element $(\overline w_1, \dots , \overline w_n)$ in $A(s_0, s_1) \otimes_k A(s_1, s_2) \otimes_k \dots \otimes_k A(s_{n-1}, s_n)$, the maps 
$$\varepsilon_n: \Hom_k(SC_n,k) \to  \Hom_{E^e}(\rad A^{\otimes n}, A)$$ given by $\varepsilon_n(f)(\overline w_1, \dots , \overline w_n)= f(s_0> s_1 > \dots > s_n)\overline{w_1 \cdots  w_n}$ commute with the corresponding boundary maps and induce the desired isomorphism of complexes.

\section{Associated incidence algebra}

\subsection{The associated incidence algebra} \label{poset}

Let $(Q,I_\nu)$ be a presentation of an algebra $A$, that is, $\nu: kQ \to A$ is a surjective morphism and $I_\nu=\Ker \nu$.  We associate to  $(Q,I_\nu)$ an incidence algebra $A(\Sigma_\nu)$, where $\Sigma_\nu$ is the poset defined as follows (see~\cite{Re2}): let $P(Q)$ be the set of paths of $Q$, let $P(Q,I_\nu)=P(Q)/\sim$ be the set of equivalence classes, where $\sim$ is the smallest equivalence relation on $P(Q)$ satisfying:
\begin{itemize}
\item[a)] If $\rho = \sum_{i=1}^m \lambda_iw_i \in I_\nu$ is a minimal relation then $w_i \sim w_j$ for all $1 \leq i,j \leq m$;
\item[b)] If $u \sim v$ then $wuw' \sim wvw'$ whenever these compositions make sense.
\end{itemize}
We denote by $[u]$ the equivalence class of a path $u$.  Observe that these are the conditions used in the definition of the fundamental group $\pi_1(Q, I_\nu)$ if we replace walks by paths. Let $\Sigma_\nu$ be the set of equivalence classes in $P(Q, I_\nu)$ that do not contain paths in $I_\nu$. Equivalent paths share source and target, and the equivalence relation is compatible with concatenation; this allows us to define $s([w])=[e_{s(w)}]$, $t([w])=[e_{t(w)}]$ and $[u][v]=[uv]$ whenever the composition makes sense and it is not equivalent to a path in $I_\nu$. The set $\Sigma_\nu$ is a poset, with $[w] \geq [w']$ if and only if there exist $[u],[v] \in \Sigma_\nu$ such that $[w']=[u][w][v]$, see~\cite{Re2}.

\begin{example}\label{ejemplo}
Consider the quiver
\[ \xymatrix{1 \ar@/^/[r]^{\alpha}\ar@/_/[r]_{\beta} & 2 \ar[r]^\gamma & 3}\]
and the ideals $I_1= < \alpha \gamma> $,  $I_2= < (\alpha - \beta) \gamma> $. The bound quivers $(Q,I_1)$ and $(Q,I_2)$ are presentations of the same algebra $A=kQ/I_1$, and the corresponding Hasse diagrams are given by
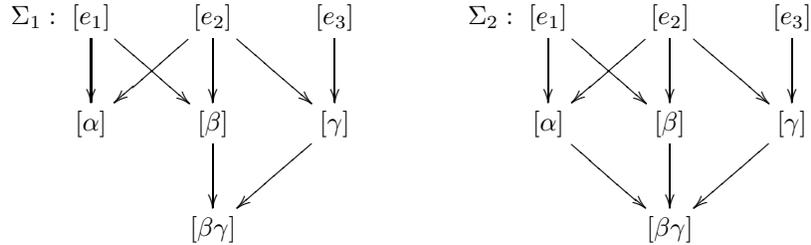
\[\Sigma_1: \xymatrix{[e_1] \ar[d] \ar[dr] & [e_2] \ar[d] \ar[dl] \ar[dr] & [e_3] \ar[d]\\
[\alpha]  & [\beta] \ar[d] & [\gamma] \ar[dl]\\
& [\beta \gamma] } \qquad \qquad \Sigma_2: \xymatrix{[e_1] \ar[d] \ar[dr] & [e_2] \ar[d] \ar[dl] \ar[dr] & [e_3] \ar[d]\\
[\alpha] \ar[dr] & [\beta] \ar[d] & [\gamma] \ar[dl]\\
& [\beta \gamma] }\]
\end{example}

\section{Homotopy coherent and compatible presentations}

In order to prove the main results in Section \ref{main} we must consider presentations satisfying two particular conditions: homotopy coherence and right (left) compatibility.  

\begin{definition}\label{P1}
A presentation $(Q,I)$ is called \textit{homotopy coherent} if for any $w, w'$ paths in $Q$ with $w \sim w'$, we have $w \in I$ if and only if $w' \in I$.
\end{definition}

This condition is necessary to construct the morphism of complexes computing $\HH^*(A(\Sigma_\nu))$ and $\HH^*(A)$, see Step 3 and Step 4 in Section \ref{4}.

\begin{remark}
If $(Q,I_\nu)$ is homotopy coherent then $\Sigma_\nu$ is just the set of equivalence classes of non-zero paths.
\end{remark}

\begin{proposition}
Let $(Q,I_\nu)$ be a presentation with $I_\nu$ generated by monomial or binomial relations.  Then $(Q,I_\nu)$ is homotopy coherent.
\end{proposition}

\begin{proof}
Let $w,w'$ be paths in $Q$ such that $w \sim w'$.  If $I_\nu$ is generated by monomial relations, then $w=w'$.  If not, there exists a finite sequence of paths
$w_0=w, w_1, \dots, w_r=w'$ such that, for all $i$ with $0 \leq i < r$, $w_i=u_i\rho_iv_i, w_{i+1}=u_i\gamma_iv_i$ with $\rho_i, \gamma_i$ appearing in a minimal relation in $I_\nu$, $u_i, v_i, \rho_i, \gamma_i$ paths in $Q$. So $\lambda_i \rho_i + \mu_i \gamma_i \in I_\nu$  for some $\lambda_i, \mu_i \in k\setminus \{0\}$.  This implies that $\lambda_i w_i + \mu_i w_{i+1} \in I_\nu$ and hence $w_i \in I_\nu$ if and only if $w_{i+1} \in I_\nu$. 
\end{proof}

Recall that an algebra $A$ is said to be \textit{schurian} if $\dim A(x,y) \leq 1$ for any $x,y \in Q_0$.

\begin{corollary}
Schurian algebras, incidence algebras and monomial algebras admit homotopy coherent presentations.
\end{corollary}

\begin{proof}
These classes of algebras admit presentations $(Q,I_\nu)$ with $I_\nu$ generated by monomial or binomial relations.
\end{proof}

\begin{definition}\label{P2}
A presentation $(Q,I_\nu)$ is called \textit{right compatible} if for any $s > s' \in SC_1(\Sigma_\nu)$ we can choose a path $u(s, s') \in Q$ such that
\begin{itemize}
\item [(i)] $s' = [v] \ s \ [u(s,s')]$; 
\item [(ii)] for any  $s > s' > s''$ in $SC_2(\Sigma_\nu)$, $u(s,s'') \sim u(s,s')u(s',s'')$.
\end{itemize} 
\end{definition}

For the sake of brevity, we refrain from stating the dual of the previous condition and leave the primal-dual translation to the reader.\\

The following result will be essentially used in Lemma \ref{left}, and shows the necessity of assuming left or right compatibility.

\begin{lemma}
If $(Q,I_\nu)$ is a right compatible presentation then there exists a family $\{ u(s, s') \in Q \vert s > s' \in SC_1(\Sigma_\nu)\}$ such that if $$[u(s_0, s_1)], \dots, [u(s_{n-1}, s_n)]$$ is the sequence associated to $s_0 > \cdots > s_n$ then $$[u(s_0, s_1)], \dots, [u(s_{i-1}, s_i)u(s_i,s_{i+1})], \dots, [u(s_{n-1}, s_n)]$$ is the sequence associated to $s_0 > \cdots > \hat{s_i} > \cdots > s_n$.
\end{lemma}

\begin{proof}
It follows by induction.
\end{proof}

\begin{example}
The presentations $(Q,I_1)$ and $(Q,I_2)$ presented in Example \ref{ejemplo} are right compatible, where the corresponding paths are given by
\begin{itemize}
\item [(i)] for $(Q,I_1)$,
\[\begin{array}{llllllll}
& u([e_1],[\alpha])=\alpha,  \quad & u([e_1],[\beta])=\beta,  \quad & u([e_2],[\alpha])=e_2,  \quad & u([e_2],[\beta])=e_2 \\ 
& u([e_2],[\gamma])=\gamma, & u([e_3],[\gamma])=e_3, & u([\beta], [\beta \gamma])=\gamma, & u([\gamma], [\beta\gamma])=e_3;
\end{array}\]
\item [(ii)] for $(Q,I_2)$,
\[ \begin{array}{llllllllll}
& u([e_1],[\alpha])=\alpha, \quad & u([e_1],[\beta])=\beta, \quad & u([e_2],[\alpha])=e_2, 
\quad & u([e_2],[\beta])=e_2 \\ 
& u([e_2],[\gamma])=\gamma, & u([e_3],[\gamma])=e_3, & u([\beta], [\beta \gamma])=\gamma, & u([\gamma], [\beta\gamma])=e_3, \\
& u([\alpha],[\beta \gamma])=\gamma.
\end{array}\]
\end{itemize} 
However the presentation $(Q,I_2)$ is not left compatible: 
\[v([e_2],[\alpha])=\alpha, \quad v([e_2],[\beta])=\beta, \quad  v([\alpha], [\beta \gamma])=e_1 \quad \mbox{and} \quad v([\beta], [\beta\gamma])=e_1.\]
On the other hand, $v([\alpha], [\beta \gamma])v([e_2],[\alpha])= \alpha$ and $v([\beta], [\beta\gamma])v([e_2],[\beta])=\beta$, hence we have no choice for $v([e_2],[\beta\gamma])$ satisfying the dual of Definition \ref{P2}(ii) since $\alpha \not \sim \beta$.
\end{example}

\begin{example} \label{ejemploNo}
Consider the quiver
\[ \xymatrix{0  \ar[r]^{\delta} & 1 \ar@/^/[r]^{\alpha}\ar@/_/[r]_{\beta} & 2 \ar[r]^\gamma & 3}\]
and the ideal $I= < \delta ( \alpha - \beta) ,  (\alpha - \beta) \gamma> $. The presentation $(Q,I)$ is neither left nor right compatible.
As in the previous example, we have no choice for $v([e_2],[\beta\gamma])$ satisfying the dual of Definition \ref{P2}(ii) and, dually, we have no choice for $u([e_1],[\delta \beta])$ satisfying Definition \ref{P2}(ii).
\end{example}

\begin{proposition}
\begin{itemize}
\item[]
\item [(a)] Monomial algebras admit right compatible presentations;
\item [(b)] Schurian algebras admit right compatible presentations;
\item [(c)] Incidence algebras admit right compatible presentations.
\end{itemize}
\end{proposition}

\begin{proof}
\begin{itemize}
\item[]
\item [(a)] The relation $\sim$ is just equality for monomial algebras;
\item [(b)] Parallel non-zero paths in a schurian algebra are linearly dependent and hence for any $s>s' \in SC_1(\Sigma_\nu)$ we have that $s' = [v] \ s \ [u]$, with $[u]$ and $[v]$ uniquely determined;
\item [(c)] Incidence algebras are schurian.
\end{itemize}
\end{proof}

\section{Main results} \label{main}

\subsection{Morphism ${\HH}^n(A(\Sigma_\nu))\to {\HH}^n(A)$}\label{4}

In order to compare the cohomological groups ${\HH}^n(A(\Sigma_\nu))$ and ${\HH}^n(A)$, we will define a morphism of complexes 
\[\Phi^*: \Hom_k(SC_*(\Sigma),k) \to \Hom_{E^e}(\rad A^{\otimes *}, A),\]
and this will be done in several steps. Recall that $F$ is the two-sided ideal of $kQ$ generated by the arrows, so $F$ is a $k$-vector space with basis the set of all non-trivial paths in $Q$. For any $n \geq 0$, let $F^{\otimes n}$ be the $n$-fold tensor product of $F$ with itself over $E$ , with ${F}^{\otimes 0}=E$.

\bigskip

\textit{Step 1}. Definition of $T_n : {F}^{\otimes n} \to SC_n(\Sigma_\nu)$.

\textit{Step 2}. Definition of $\partial_{n+1} : {F}^{\otimes n+1} \to  {F}^{\otimes n}$  and $\delta_{n+1} : SC_{n+1}(\Sigma_\nu) \to SC_n(\Sigma_\nu)$.

\textit{Step 3}. Proof of the commutativity of the diagram
\[ \xymatrix{ {F}^{\otimes n+1} \ar[rr]^{T_{n+1}} \ar[d]^{\partial_{n+1}} & & SC_{n+1}(\Sigma_\nu) \ar[d]^{\delta_{n+1}} \\
{F}^{\otimes n} \ar[rr]^{T_n}  & & SC_n(\Sigma_\nu) } \]

\textit{Step 4}. Description of the morphism $\Phi^n : \Hom_k(SC_n,k) \to \Hom_{E^e} (\rad A^{\otimes n}, A)$.

\bigskip

\noindent \textit{Step 1}. Let $T_n : {F}^{\otimes n} \to SC_n(\Sigma_\nu)$ be the $k$-linear map defined inductively by:
\begin{itemize}
\item [(i)] $T_0(e_i)=[e_i]$;
\item [(ii)] For any path $w$ in $F$, put 
\[T_1(w)  =  [e_{s(w)}]>[w]  - [e_{t(w)}]>[w] \]
if $[w] \in \Sigma_\nu$, and zero otherwise;
\item [(iii)] For any basis element $(w_1, \dots, w_n)$ in ${F}^{\otimes n}$, put 
\begin{multline*}
T_n(w_1, \dots, w_n) \\
 =  \left[ T_{n-1}(w_1, \dots, w_{n-1}) + (-1)^n T_{n-1}(w_2, \dots, w_n) \right]  > [w_1 \cdots w_n]
\end{multline*}
if  $[w_1 \dots w_{n}] \in \Sigma_\nu$, and zero otherwise.
\end{itemize}
Note that we put
$\left[ \sum_{i} \lambda_i \ s_0^i > \cdots > s_{n-1}^i \right] > s_n :=  \sum_{i} \lambda_i \ s_0^i > \cdots > s_{n-1}^i > s_n$.

\medskip

\noindent \textit{Step 2}. For $n \geq 0$, let $\partial_{n+1} : {F}^{\otimes n+1} \to  {F}^{\otimes n}$ be the $k$-linear map defined by
\begin{equation*}
\begin{split}
\partial_1 (w) & =e_{t(w)}-e_{s(w)}, \\
\partial_{n+1} (w_0, \dots , w_n) & =  (w_1, \dots , w_n) + \widetilde{\partial_{n+1}} (w_0, \dots , w_n) + (-1)^{n+1} (w_0, \dots , w_{n-1}) \\
& =  (w_1, \dots , w_n) + \sum_{i=1}^{n} (-1)^i (w_0, \dots, w_{i-1}w_i, \dots , w_n)  \\
& \quad +  (-1)^{n+1} (w_0, \dots , w_{n-1}),
\end{split}
\end{equation*}
and let $\delta_{n+1} : SC_{n+1}(\Sigma_\nu) \to SC_n(\Sigma_\nu)$ be the $k$-linear map defined by
\[\delta_{n+1} (s_0 > \cdots > s_{n+1})  = \sum_{i=0}^{n+1} (-1)^i s_0 > \cdots > \hat{s_i} > \cdots > s_{n+1}.\]

\medskip

\noindent \textit{Step 3}.
From now on, we will consider homotopy coherent presentations $(Q,I)$.
We stress that, under this assumption, $T_n$ is compatible with the equivalence relation used to define $\Sigma_\nu$, that is, if $w_i \sim u_i$ for all $i$ with $1 \leq i \leq n$ then $T_n(w_1, \dots, w_n) = T_n(u_1, \dots, u_n)$.

\begin{remark}\label{trivial}
If $(Q,I_\nu)$ is homotopy coherent and $[u],[v] \in \Sigma_\nu$ are such that $[u]=[v][u]$ then $v$ is a trivial path. In fact, if $v$ is a path of positive length, $u \sim vu$ implies that $u \sim v^m u$ and, for $m$ sufficiently large, the path $v^m u$ belongs to the admissible ideal $I_\nu$, a contradiction.
\end{remark}

\begin{lemma}\label{h}
Let $(Q,I_\nu)$ be a homotopy coherent presentation of $A$.  Then \\
$T_{n}  \partial_{n+1} (w_0, \dots , w_n)=  \delta_{n+1} T_{n+1}(w_0, \dots , w_n)$, for any $n \geq 0$ and $(w_0, \dots , w_n)$ a basis element in $F^{\otimes n+1}$ with $w_0 \cdots w_n \not \in I_\nu$.
\end{lemma}

\begin{proof}
Observe that $w_0 \cdots w_n \not \in I_\nu$ implies that 
$w_1 \cdots w_n \not \in I_\nu$ and $w_0 \cdots w_{n-1} \not \in I_\nu$,
and the homotopy coherence of the presentation implies that $[w_0 \dots w_{n}] \in \Sigma_\nu$ if and only if $w_0 \cdots w_n \not \in I_\nu$. 
A direct computation shows that
$T_0  \partial_1 (w) =  [e_{t(w)}]-[e_{s(w)}] =  \delta_1 T_1 (w)$ and
\begin{equation*}
\begin{split}
T_1  \partial_2 (w_0, w_1) & = T_1(w_1)-[e_{s(w_0)}]>[w_0w_1] + [e_{t(w_1)}]>[w_0w_1] +T_1(w_0)  \\
& =  \delta_2 T_2(w_0, w_1)
\end{split}
\end{equation*}
for any $w \not \in I_\nu$, $w_0w_1 \not \in I_\nu$.
Now we proceed by induction, assuming that the desired equality holds for any $j$ such that $0 \leq j < n$, $n>1$. 
\begin{equation*}\begin{split}
& T_n   \partial_{n+1}  (w_0,  \dots , w_n)\\
& \quad =  T_n(w_1, \dots, w_n) - T_n(w_0w_1, \dots, w_n) + (-1)^n T_n(w_0, \dots, w_{n-1}w_n) \\
& \qquad  + (-1)^{n+1}T_n(w_0, \dots, w_{n-1}) - T_n (w_0,\widetilde{\partial_{n-1}} (w_1, \dots, w_{n-1}),w_n) \\
& \quad =  T_n(w_1, \dots, w_n) - T_n(w_0w_1, \dots, w_n)  + (-1)^n T_n(w_0, \dots, w_{n-1}w_n)  \\
& \qquad +  (-1)^{n+1}T_n(w_0, \dots, w_{n-1}) \\
& \qquad  - T_{n-1}(w_0, \widetilde{\partial_{n-1}} (w_1, \dots, w_{n-1})) > [w_0 \cdots w_n] \\
& \qquad - (-1)^{n} T_{n-1}(\widetilde{\partial_{n-1}} (w_1, \dots, w_{n-1}),w_n) > [w_0 \cdots w_n] \\
& \quad =  T_n(w_1, \dots, w_n) - T_n(w_0w_1, \dots, w_n)  + (-1)^n T_n(w_0, \dots, w_{n-1}w_n) \\
& \qquad +(-1)^{n+1}T_n(w_0, \dots, w_{n-1}) \\
& \qquad  +  \left[ T_{n-1} \partial_{n} (w_0, \dots, w_{n-1})  
+ (-1)^{n+1}  T_{n-1} \partial_{n} (w_1, \dots , w_{n-1},w_n) \right] > [w_0 \cdots w_n]  \\
& \qquad  +  \left[-T_{n-1}(w_1, \dots, w_{n-1})  + (-1)^n T_{n-1} (w_2, \cdots, w_n)\right] > [w_0 \cdots w_n] \\
& \qquad  +  \left[T_{n-1}(w_0w_1, \dots, w_{n-1})  - T_{n-1} (w_1, \dots , w_{n-1}w_n)\right] > [w_0 \cdots w_n] \\
& \qquad  +  \left[ (-1)^{n+1} T_{n-1}(w_0, \dots , w_{n-2})  + T_{n-1} (w_1, \dots , w_{n-1})\right] > [w_0 \cdots w_n].
\end{split}
\end{equation*}
By the inductive hypothesis and using the inductive definition of $T_n$ we have 
\begin{equation*}
\begin{split}
& T_n  \partial_{n+1} (w_0, \dots , w_n) \\
& \quad = T_n(w_1, \dots , w_n) +(-1)^{n+1}T_n(w_0, \dots , w_{n-1}) \\
 & \qquad + \left[  \delta_{n}T_{n} (w_0, \dots , w_{n-1}) +  (-1)^{n+1}    \delta_{n}T_{n} (w_1, \dots , w_n) \right] > [w_0 \cdots w_n] \\
& \quad =  \delta_{n+1} \left( T_n(w_0, \dots , w_{n-1})> [w_0 \cdots w_n]\right)\\
& \qquad +  (-1)^{n+1}  \delta_{n+1} \left(T_n (w_1, \dots , w_n) > [w_0 \cdots w_n] \right)\\
& \quad =  \delta_{n+1} T_{n+1}(w_0, \dots , w_n).
\end{split} 
\end{equation*}
\end{proof}

\bigskip

\noindent \textit{Step 4}. We are now in a position to describe the morphisms 
\[\Phi^n : \Hom_k(SC_n(\Sigma_\nu),k) \to \Hom_{E^e} (\rad A^{\otimes n}, A).\]
Note that $\rad A \simeq F/I_\nu$ and $\rad A^{\otimes n} \simeq F^{\otimes n}/R$ with 
\[R= \sum_{s=0}^{n-1} F^{\otimes s} \otimes_E I_\nu \otimes_E F^{\otimes n-s-1}.\]
For any $w \in F$ we write $\overline w$ for the class of $w$ modulo $I_\nu$.  For any $f$ in $\Hom_k(SC_n(\Sigma_\nu),k)$ let $\widetilde \Phi^n (f) : F^{\otimes n} \to A$ be the $k$-linear map defined by
\[\widetilde \Phi^n(f) (w_1, \dots , w_n) = f(T_n(w_1, \dots , w_n)) \overline{w_1 \cdots w_n},\]
where $(w_1, \dots , w_n)$ is a basis element in $F^{\otimes n}$.

\medskip

Let $(w_1, \dots , \rho, \dots , w_n) \in F^{\otimes s} \otimes_E I_\nu \otimes F^{\otimes n-s-1}$. If $\rho$ is a path, it is clear that 
$\widetilde \Phi^n(f)(w_1, \dots , \rho, \dots , w_n)=0$.  If $\rho= \sum_{i=1}^m \lambda_i v_i$ is a minimal relation, $m > 1$, then 
\[T_n (w_1, \dots , v_i, \dots , w_n) = T_n (w_1, \dots , v_1, \dots , w_n)\]
for any $i$ with $1 \leq i \leq m$. Hence
\begin{equation*}
\begin{split}
\widetilde \Phi^n(f)(w_1, \dots , \rho, \dots , w_n) & =  \sum_{i=1}^m \lambda_i \widetilde \Phi^n(f)(w_1, \dots , v_i, \dots, w_n) \\
&=  \sum_{i=1}^m \lambda_i f(T_n(w_1, \dots, v_i, \dots, w_n)) \overline{w_1 \cdots v_i \cdots  w_n} \\
&=  f(T_n(w_1, \dots, v_1, \dots, w_n)) \sum_{i=1}^m \lambda_i \overline{w_1 \cdots v_i \cdots  w_n} = 0.
\end{split}
\end{equation*}
The ideal $I_\nu$ is generated by paths and minimal relations, so $\widetilde \Phi^n(f)(R)=0$, and then $\widetilde \Phi^n$ induces a map $\Phi^n(f): \rad A^{\otimes n} \to A$ given by 

\begin{equation*}
\begin{split}
\Phi^0(f)(e_i) & = f(T_0(e_i))e_i= f([e_i])e_i, \\
\Phi^n(f) (\overline w_1, \dots , \overline w_n)& = f(T_n(w_1, \dots , w_n)) \overline{w_1 \cdots w_n}.
\end{split}
\end{equation*}

\begin{proposition} If $(Q,I_\nu)$ is a homotopy coherent presentation of $A$, the map $\Phi^* : \Hom_k(SC_*(\Sigma_\nu),k) \to \Hom_{E^e} (\rad A^{\otimes *}, A)$ is a morphism of complexes and, hence, it induces morphisms ${\HH}(\Phi^n): {\HH}^n(A(\Sigma_\nu)) \to {\HH}^n(A)$, for any $n \geq 0$.
\end{proposition}

\begin{proof}
We have to show that, for any $n \geq 0$, the diagram
\[\xymatrix{ \Hom_k(SC_n(\Sigma_\nu),k) \ar[r]^{\Phi^n} \ar[d]^{B^n}  &  \Hom_{E^e} (\rad A^{\otimes n}, A) \ar[d]^{b^n} \\ \Hom_k(SC_{n+1}(\Sigma_\nu),k) \ar[r]^{\Phi^{n+1}} & \Hom_{E^e} (\rad A^{\otimes n+1}, A) } \]
is commutative.  A direct computation shows that
\begin{equation*}
\begin{split}
\Phi^{n+1}(B^{n}f) (\overline w_0, \dots , \overline w_n) & =  (B^{n}f)(T_{n+1}(w_0, \dots , w_n)) \overline{w_0 \cdots w_n}\\
& =  f(\delta_{n+1} T_{n+1}(w_0, \dots , w_n)) \overline{w_0 \cdots w_n}
\end{split}
\end{equation*}
and 
\[(b^n \Phi^n(f)) (\overline w_0, \dots , \overline w_n)  = f(T_n \partial_{n+1}(w_0, \dots , w_n)) \overline{w_0 \cdots w_n}.\]
The desired equality follows from Lemma~\ref{h}.
\end{proof}

\subsection{The complex $\Ker \Phi^*$}

We will show that the complex $\Ker \Phi^*$ is exact, by constructing a contraction homotopy $S_n:  \Ker \Phi^n \to \Ker \Phi^{n-1}$. Observe that 
\begin{equation*}
\begin{split}
\Ker \Phi^n & =  \{ f \in \Hom_k(SC_n(\Sigma_\nu),k) : f(T_n(w_1, \dots , w_n)) \overline{w_1 \cdots w_n} = 0,  \\
& \qquad \quad \mbox{for any basis element $(w_1, \dots , w_n) \in F^{\otimes n}$}\} \\
& =  \{ f \in \Hom_k(SC_n(\Sigma_\nu),k) : f(T_n(w_1, \dots , w_n)) = 0,  \\
& \qquad \quad \mbox{for any basis element $(w_1, \dots , w_n) \in F^{\otimes n} \setminus R $}\}.
\end{split}
\end{equation*}
In order to construct the homotopy, the chosen presentation $(Q,I_\nu)$ of $A$ must be left or right compatible, see Definition \ref{P2}. 
In this case we choose a family 
\[\{ u(s, s') \in Q \vert s > s' \in SC_1(\Sigma_\nu)\}\]
satisfying Definition \ref{P2} and let $G^0 : SC_0(\Sigma_\nu) \to SC_1(\Sigma_\nu)$ be the map given by 
\[ G^0([w])=
\begin{cases}
0 & \mbox{if $w \in E$}, \\
[e_{t(w)}] > [w] & \mbox{if $w \in F$}.
\end{cases} \]
Let $G^{n}: SC_{n}(\Sigma_\nu) \to SC_{n+1}(\Sigma_\nu)$, for $n > 0$, be defined inductively in the following way: for any $s_0 > \cdots > s_n$ in $SC_n(\Sigma_\nu)$, denote 
\[u_i = u(s_{i-1}, s_{i}), \qquad \Omega= \Omega(s_0 > \cdots > s_n)=\{ i : u_i \in E , 0< i \leq n\},\]
and put
\begin{itemize}
\item [(i)] $G^{n} (s_0 > \cdots > s_n) = G^{n-1}(s_0 > \cdots > s_{n-1})> s_n$ if $\Omega \not = \emptyset$ or 
$[u_1 \cdots u_n] = s_n$;
\item [(ii)] $G^{n} (s_0 > \cdots > s_n) = \left[ G^{n-1}(s_0 > \cdots > s_{n-1}) +(-1)^{n} T_{n}(u_1, \dots, u_n) \right] > s_n$ otherwise.
\end{itemize}

\medskip

In order to prove that the complex $\Ker \Phi^*$ admits a contraction homotopy, we need the following lemma.

\begin{lemma} \label{left} Let $(Q, I_\nu)$ be a homotopy coherent, right compatible presentation of $A$. Then
\begin{itemize}
\item [(i)] $ \delta_1 G^0 ([w])= [w]- T_0(e_{t(w)})$, for any $[w] \in SC_0(\Sigma_\nu)$;
\item [(ii)] $(\delta_{n+1} G^{n} + G^{n-1}  \delta_n )(s_0 > \cdots > s_n)= s_0 > \cdots > s_n$ if $\Omega(s_0 > \cdots > s_n) \not = \emptyset$;
\item [(iii)] $(\delta_{n+1} G^{n} + G^{n-1} \delta_n )(s_0 > \cdots > s_n)= s_0 > \cdots > s_n - T_{n}(u_1, \dots, u_n)$, otherwise.
\end{itemize}
\end{lemma}

\begin{proof}
A direct computation shows that $\delta_1 G^0 ([w])= [w]- T_0(e_{t(w)})$
and that the assertion is true for $\delta_2 G^1 + G^0  \delta_1$. Now assume that $n > 1$ and proceed by induction.
Let $s_0 > \cdots > s_n \in SC_n(\Sigma_\nu)$ and consider the following cases:
\begin{itemize}
\item [ a)] $\Omega(s_0 > \cdots > s_n) = \emptyset$, $[u_1 \cdots u_n] \not = s_n$;
\item [ b)] $\Omega(s_0 > \cdots > s_n) = \emptyset$, $[u_1 \cdots u_n] = s_n$;
\item [ c)] \begin{itemize} \item [1)] $\Omega(s_0 > \cdots > s_n)  = \{ n \}$;
\item [2)] $\Omega(s_0 > \cdots > s_n)  = \{ i \}$, $0< i < n$;
\item [3)] $\Omega(s_0 > \cdots > s_n)$ has at least two elements.
\end{itemize}
\end{itemize}
\textit{Case} (a). If $s_0 > \cdots > s_n$  is such that $\Omega(s_0 > \cdots > s_n) = \emptyset$ and $[u_1 \cdots u_n] \not = s_n$, we have
\begin{equation*} \begin{split}
& \delta_{n+1} G^{n} (s_0 > \cdots > s_n) \\
& \quad =   \delta_{n+1} (G^{n-1}(s_0 > \cdots > s_{n-1}) > s_n)  +  (-1)^{n}  \delta_{n+1} (T_{n}(u_1, \dots, u_n) > s_n)  \\
& \quad =  \delta_{n} G^{n-1}(s_0 > \cdots >  s_{n-1})  > s_n  + (-1)^{n+1}  G^{n-1}(s_0 > \cdots > s_{n-1}) \\
& \qquad +   (-1)^{n}   \delta_{n} T_{n}(u_1, \dots, u_n) > s_n  +   (-1)^n (-1)^{n+1} T_{n}(u_1, \dots, u_n).
\end{split}
\end{equation*}
Using the inductive hypothesis and Lemma~\ref{h} we get
\begin{equation*} \begin{split}
& \delta_{n+1} G^{n} (s_0 > \cdots > s_n) \\
& \quad=  s_0 > \cdots > s_n - G^{n-2} \delta_{n-1} (s_0 > \cdots >  s_{n-1}) > s_n \\
& \qquad -  T_{n-1}(u_1 , \dots, u_{n-1}) > s_n  +  (-1)^{n+1}  G^{n-1}(s_0 > \cdots > s_{n-1})  \\
& \qquad +   (-1)^{n}  T_{n-1} \partial_{n}(u_1, \dots, u_n) > s_n  -   T_{n}(u_1, \dots, u_n) \\
& \quad=  s_0 > \cdots > s_n  - T_{n}(u_1, \dots, u_n)  -  G^{n-1} \delta_n(s_0 > \cdots > s_n).
\end{split}
\end{equation*}
\textit{Case} (b). If $s_0 > \cdots > s_n$  is such that $\Omega(s_0 > \cdots > s_n) = \emptyset$ and $[u_1 \cdots u_n] = s_n$ then
\begin{equation*} \begin{split}
& \delta_{n+1} G^{n} (s_0 > \cdots > s_n)  \\
& \quad =   \delta_{n+1} (G^{n-1}(s_0 > \cdots > s_{n-1}) > s_n) \\
& \quad =   \delta_{n} G^{n-1}(s_0 > \cdots >  s_{n-1})  > s_n  +  (-1)^{n+1}  G^{n-1}(s_0 > \cdots > s_{n-1}) \\
& \quad =  s_0 > \cdots > s_n - G^{n-2} \delta_{n-1} (s_0 > \cdots >  s_{n-1}) > s_n \\
& \qquad -  T_{n-1}(u_1 , \dots, u_{n-1}) > s_n  +  (-1)^{n+1}  G^{n-1}(s_0 > \cdots > s_{n-1}).
\end{split}
\end{equation*}
On the other hand
\begin{equation*} \begin{split}
& G^{n-1} \delta_{n}(s_0 > \cdots > s_n) \\
& \qquad =  G^{n-1}(\delta_{n-1} (s_0 > \cdots > s_{n-1}) > s_n)  +  (-1)^{n}  G^{n-1}(s_0 > \cdots > s_{n-1}) \\
& \qquad =  G^{n-2}\delta_{n-1} (s_0 > \cdots > s_{n-1}) > s_n + (-1)^{n-1} T_{n-1}(u_2, \dots, u_n) > s_n \\
& \qquad \quad +  (-1)^{n}  G^{n-1}(s_0 > \cdots > s_{n-1}). \end{split}
\end{equation*}
Hence
\begin{equation*} \begin{split}
&(\delta_{n+1} G^{n}+ G^{n-1} \delta_{n})(s_0 > \cdots > s_n) \\
& \quad =  s_0 > \cdots > s_n  -  \left[ T_{n-1}(u_1 , \dots, u_{n-1}) + (-1)^n T_{n-1}(u_2, \dots, u_n)\right] > s_n \\
& \quad =   s_0 > \cdots > s_n - T_n(u_1, \dots , u_n). \end{split}
\end{equation*}
\textit{Case} (c). If $s_0 > \cdots > s_n$ is such that $\Omega(s_0 > \cdots > s_n) \not = \emptyset$ then 
\begin{equation*} \begin{split}
& \delta_{n+1} G^{n} (s_0 > \cdots > s_n) \\
& \quad =  \delta_{n+1} (G^{n-1}(s_0 > \cdots > s_{n-1}) > s_n) \\
& \quad =  \delta_{n} G^{n-1}(s_0 > \cdots >  s_{n-1})  > s_n  +  (-1)^{n+1}  G^{n-1}(s_0 > \cdots > s_{n-1}). \end{split}
\end{equation*}
\textit{Case} (c 1). If $\Omega(s_0 > \cdots > s_n)  = \{ n \}$ then $u(s_{n-2},s_n)=u_{n-1}$ and
\[\Omega(s_0 > \cdots> \hat s_j > \cdots > s_n) = \begin{cases} \emptyset \quad & \mbox{if $j=n-1,n$}, \\ \{ n \} & \mbox{otherwise}. \end{cases}\]
Moreover $[u_1 \cdots u_{n-1}] \not = s_n$.  In fact, if $[u_1 \cdots u_{n-1}] = s_n$, using Remark~\ref{trivial} we deduce that $s_n=s_{n-1}$, a contradiction.  So
\begin{equation*} \begin{split}
& G^{n-1} \delta_{n}(s_0 > \cdots > s_n) \\
& \quad =  G^{n-1}(\delta_{n-1} (s_0 > \cdots > s_{n-1}) > s_n) +  (-1)^{n}  G^{n-1}(s_0 > \cdots > s_{n-1}) \\
& \quad =  G^{n-2}\delta_{n-1} (s_0 > \cdots > s_{n-1}) > s_n  + T_{n-1}(u_1, \dots, u_{n-1}) > s_n \\
& \qquad +  (-1)^{n}  G^{n-1}(s_0 > \cdots > s_{n-1}). \end{split}
\end{equation*}
\textit{Case} (c 2). If $\Omega(s_0 > \cdots > s_n)  = \{ i \}$, $0< i < n$, then 
\[\Omega(s_0 > \cdots> \hat s_j > \cdots > s_n) = \begin{cases} \emptyset \quad & \mbox{if $j=i-1,i$}, \\ \{ i \} & \mbox{otherwise}, \end{cases}\]
\begin{multline*}
 G^{n-1}(s_0> \cdots > \hat s_{i-1} > \cdots > s_n) - G^{n-1}(s_0> \cdots > \hat s_i > \cdots > s_n)    \\
 = \left[ G^{n-2}(s_0> \cdots > \hat s_{i-1} > \cdots > s_{n-1}) - G^{n-2}(s_0> \cdots > \hat s_i > \cdots > s_{n-1})\right] >s_n
\end{multline*}
and so
\begin{equation*} \begin{split}
& G^{n-1} \delta_{n}(s_0 > \cdots > s_n)  \\
& \quad =  G^{n-1}(\delta_{n-1} (s_0 > \cdots > s_{n-1}) > s_n)  +  (-1)^{n}  G^{n-1}(s_0 > \cdots > s_{n-1}) \\
& \quad =  G^{n-2}\delta_{n-1} (s_0 > \cdots > s_{n-1}) > s_n  + (-1)^{n}  G^{n-1}(s_0 > \cdots > s_{n-1}). \end{split}
\end{equation*}
\textit{Case} (c 3). If $\Omega(s_0 > \cdots > s_n)$ has at least two elements then 
\[\Omega(s_0 > \cdots > \hat s_i > \cdots > s_n) \not = \emptyset\] and 
\begin{equation*} \begin{split}
& G^{n-1} \delta_{n}(s_0 > \cdots > s_n) \\
& \quad =  G^{n-1}(\delta_{n-1} (s_0 > \cdots > s_{n-1}) > s_n)  +  (-1)^{n}  G^{n-1}(s_0 > \cdots > s_{n-1}) \\
& \quad =  G^{n-2}\delta_{n-1} (s_0 > \cdots > s_{n-1}) > s_n  + (-1)^{n}  G^{n-1}(s_0 > \cdots > s_{n-1}). \end{split}
\end{equation*}
Then, in cases (c 2) and (c 3),
\begin{multline*}
(\delta_{n+1} G^{n}+G^{n-1} \delta_{n}) (s_0 > \cdots > s_n) \\
= \left[ (\delta_{n} G^{n-1}+G^{n-2}\delta_{n-1})(s_0 > \cdots >  s_{n-1}) \right] > s_n 
\end{multline*}
and in case (c 1)
\begin{multline*}
 (\delta_{n+1} G^{n}+G^{n-1} \delta_{n}) (s_0 > \cdots > s_n) \\
 =  \left[ (\delta_{n} G^{n-1}+G^{n-2}\delta_{n-1})(s_0 > \cdots >  s_{n-1}) + T_{n-1}(u_1, \dots, u_{n-1}) \right] > s_n. 
\end{multline*}
Hence the assertion follows by induction.
\end{proof}

Note that $G^{n}T_n =0$: a direct computation proves this for $n=0,1$, and an inductive procedure completes the proof since
\[T_n(w_1, \dots, w_n) = [e_{s(w_1)}]>[w_1]> \cdots >[w_1 \cdots w_n] + \mbox{simplices with $\Omega \not = \emptyset$}\]
and hence
\begin{equation*} \begin{split}
& G^n T_n (w_1, \dots, w_n) \\
& \quad = G^n ( \left[ T_{n-1} (w_1, \dots, w_{n-1})  + (-1)^n T_{n-1} (w_2, \dots, w_n) \right] > [w_1 \cdots w_n])\\
& \quad =  G^{n-1} ( T_{n-1} (w_1, \dots, w_{n-1})  +(-1)^n  T_{n-1} (w_2, \dots, w_n)) > [w_1 \cdots w_n].
\end{split}
\end{equation*}
So, for any $f \in \Ker \Phi^{n+1}$, the composition $f G^{n}$ belongs to $\Ker \Phi^{n}$. 
Let 
\[S_{n+1}: \Ker \Phi^{n+1} \to \Ker \Phi^{n}\] be the map defined by $S_{n+1}(f)= f G^{n}$.

\begin{proposition} Let $(Q, I_\nu)$ be a homotopy coherent, right compatible presentation of $A$. The map $S_{n+1}: \Ker \Phi^{n+1} \to \Ker \Phi^{n}$ is a contraction homotopy.
\end{proposition}

\begin{proof}
It follows immediately from the previous lemma since
\begin{equation*}
\begin{split}
(S_1B^0)(f) & = B^0 f G^ 0 = f \delta_1 G^0=f, \quad \mbox{and} \\
(S_{n+1} B^{n} + B^{n-1} S_{n})(f)&= B^{n} f  G^{n} + B^{n-1} f G^{n-1}=
f \delta_{n+1} G^n + f G^{n-1} \delta_{n}=f
\end{split}
\end{equation*}
because $f T_{0}=0=f T_{n}$.
\end{proof}

\begin{theorem}\label{teorema}
Let $(Q,I_\nu)$ be a homotopy coherent, right compatible presentation of an algebra $A$.  If $\Phi^{n-1}$ is a surjective morphism, then ${\HH}(\Phi^n): {\HH}^n(A(\Sigma_\nu)) \to {\HH}^n(A)$ is an injective morphism.
\end{theorem}

\begin{proof}
It follows by diagram chasing of elements in the commutative diagram of complexes computing the corresponding cohomology groups.
\end{proof}

Monomial algebras without non-zero oriented cycles, schurian algebras and incidence algebras satisfy the assumptions of the following corollary.

\begin{corollary}\label{46}
Let $(Q,I_\nu)$ be a homotopy coherent, right compatible presentation of an algebra $A$. If $\dim_k A(x,x)=1$ for any $x \in Q_0$ then the morphism
${\HH}(\Phi^1): {\HH}^1(A(\Sigma_\nu)) \to {\HH}^1(A)$ is injective.
\end{corollary}

\begin{proof}
The proof follows from the previous theorem by observing that the assumed hypotheses imply that $\Phi^0$ is a surjective map.
\end{proof}

\begin{corollary} \label{caso1}
If $A$ is an incidence algebra, then ${\HH}(\Phi^n): {\HH}^n(A(\Sigma_\nu)) \to {\HH}^n(A)$ is an isomorphism for any $n \geq 0$.
\end{corollary}

\begin{proof}
The assertion is clear for $n=0$. Since $A={\Ia}(\Sigma)$ is an incidence algebra, recall from Section \ref{posetincidencia} that 
\[\rad A^{\otimes n}=  \oplus_{s_0> s_1 > \dots > s_n} A(s_0, s_1) \otimes_k A(s_1, s_2) \otimes_k \dots \otimes_k A(s_{n-1}, s_n),\]
with $\dim_k A(s_{i-1}, s_i)=1$ for all $i$ with $1\leq i \leq n$. Taking a set of basis elements $(\overline w_1, \dots , \overline w_n)$ in $A(s_0, s_1) \otimes_k A(s_1, s_2) \otimes_k \dots \otimes_k A(s_{n-1}, s_n)$, the maps $g: \rad A^{\otimes n} \to A$ defined by
\[g(\overline v_1, \dots , \overline v_n)= \begin{cases}
\lambda \overline{ w_1 \cdots w_n} \quad & \mbox{if $(\overline v_1, \dots , \overline v_n)=(\overline w_1, \dots , \overline w_n)$},\\
0 & \mbox{otherwise},
\end{cases}\]
form a basis of the $k$-vector space $\Hom_{E^e} (\rad A^{\otimes n}, A)$.  Take $f \in \Hom_k(SC_n(\Sigma_\nu),k)$ defined as follows:
\[f(s_0 > \cdots > s_n)= \begin{cases}
\lambda \quad & \mbox{if $s_0 > \cdots > s_n = [e_{s(w_1)}]>[w_1]>\cdots > [w_1 \cdots w_n]$},\\
0 & \mbox{otherwise}.
\end{cases}\]
Now $\Phi^n(f)=g$ and hence $\Phi^{n}$ is a surjective map for any $n \geq 0$. Then we get the short exact sequence of the complexes
\[0 \to \Ker (\Phi^*) \to \Hom_k(SC_*(\Sigma_\nu),k) \to \Hom_{E^e}(\rad A^{\otimes *},A) \to 0,\]
so that we get the long exact sequence
\[\dots \to {\HH}^n(\Ker \Phi^*) \to {\HH}^n(A(\Sigma_\nu)) \to {\HH}^n(A) \to {\HH}^{n+1}(\Ker \Phi^*) \to \dots\]
\end{proof}

\section{Examples}

\begin{example}
Consider the presentations $(Q,I_1)$ and $(Q,I_2)$ of the algebra $A=kQ/I_1$ given by
\[ Q:  \xymatrix{1 \ar@/^/[r]^{\alpha}\ar@/_/[r]_{\beta} & 2 \ar[r]^\gamma & 3, }\]
$I_1= < \alpha \gamma> $ and $I_2= < (\alpha - \beta) \gamma> $, presented in Example~\ref{ejemplo}. Using \cite{IZ} we construct the reduced posets $\overline \Sigma_1, \overline \Sigma_2$, as described in the Introduction, 
\[\overline \Sigma_1: \xymatrix{[e_1] \ar[d] \ar[dr] & [e_2] \ar[d] \ar[dl] & [e_3] \ar[dl]\\
[\alpha]  & [\beta\gamma]  & } \qquad \qquad \overline \Sigma_2: \xymatrix{[e_1] \ar[dr] & [e_2] \ar[d] & [e_3] \ar[dl]\\
& [\beta \gamma] }\]
and from \cite[2.5,2.2]{GR} and \cite[5.3]{H} we get
\[{\HH}^i(A(\Sigma_1)) =
\begin{cases}
k \quad &\mbox{if $i=0, 1$}, \\
0 \quad &\mbox{otherwise}
\end{cases}, \qquad
{\HH}^i(A(\Sigma_2)) =
\begin{cases}
k \quad &\mbox{if $i=0$}, \\
0 \quad &\mbox{otherwise}.
\end{cases}\]
From \cite[5.3,1.6]{H} we get
\[{\HH}^i(A) =
\begin{cases}
k \quad &\mbox{if $i=0$}, \\
k^2 \quad &\mbox{if $i=1$}, \\
0 \quad &\mbox{otherwise}.
\end{cases} \]
\end{example}

Now we present families of algebras where the non-vanishing of some Hochschild cohomology groups can be deduced.

\begin{example}
Consider the quiver
\[ \xymatrix{1 \ar @/^/ @{-} [rr]^{\alpha_1}\ar @/_/ @{-} [rr]_{\alpha_n} & \vdots & 2,}\]
with arrows $\alpha_1, \dots \alpha_n$ with any orientation, and let $A=kQ/F^2$.  The quiver of the corresponding incidence algebra $A(\Sigma)$ is given by
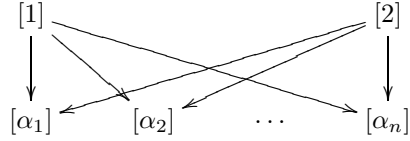
\[\xymatrix{ [1] \ar[d] \ar[rd] \ar[drrr] & & & [2] \ar[d] \ar[dll] \ar[dlll] \\
[\alpha_1] & [\alpha_2] &  \dots & [\alpha_n]} \]
and from~\cite[1.6]{H} we get 
\[{\HH}^i(A(\Sigma)) =
\begin{cases}
k \quad &\mbox{if $i=0$}, \\
k^{n-1} \quad &\mbox{if $i=1$}, \\
0 \quad &\mbox{otherwise}.
\end{cases}\]
Hence $\dim_k {\HH}^1(A) \geq n-1$ by Corollary \ref{46}.
If $\alpha_1, \cdots, \alpha_n$ share starting and ending points, then from~\cite[1.6]{H} we get 
\[{\HH}^i(A) =
\begin{cases}
k \quad &\mbox{if $i=0$}; \\
k^{n^2-1} \quad &\mbox{if $i=1$}; \\
0\quad &\mbox{otherwise}.
\end{cases}\]
The particular case $n=2$, $\alpha_1,\alpha_2$ with opposite orientations, has been considered in~\cite{C2} and 
\[{\HH}^i(A) =
\begin{cases}
k \quad &\mbox{if $i=0, 4s, 4s+1$}; \\
0\quad &\mbox{otherwise}.
\end{cases} \]
\end{example}

The previous example shows that the injective morphism described in Theorem~\ref{teorema} can not be expected to be an isomorphism in general. It would be nice to determined all the algebras that admit a distinguished presentation making the mentioned morphism an isomorphism. 

\begin{example}
Let $Q_n$ be the quiver
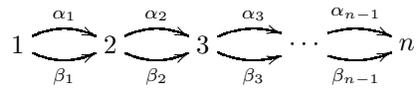
\[\xymatrix{ 1 \ar @/^/ [r]^{\alpha_1}\ar @/_/ [r]_{\beta_1} & 2 \ar @/^/ [r]^{\alpha_2}\ar @/_/ [r]_{\beta_2} &  3  \ar @/^/ [r]^{\alpha_3}\ar @/_/  [r]_{\beta_3} &    \cdots  \ar @/^/ [r]^{\alpha_{n-1}}\ar @/_/ [r]_{\beta_{n-1}}  & n}\]
and let $A_n=kQ_n/F^2$.  The quiver of the corresponding incidence algebra $A(\Sigma_n)$ is given by
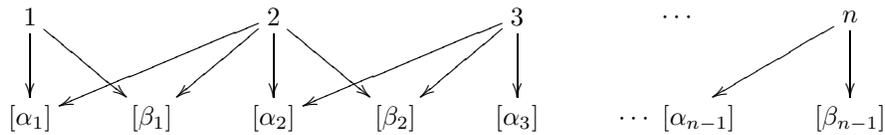
\[\xymatrix{ 1 \ar[d] \ar[rd] & & 2 \ar[dll] \ar[dl] \ar[d] \ar[dr]  & & 3 \ar[dll] \ar[dl] \ar[d]  &  \cdots &  n \ar[dl] \ar[d]\\
[\alpha_1] &  [\beta_1] &  [\alpha_2] & [\beta_2] &   [\alpha_{3}] & \cdots  \  [\alpha_{n-1}] & [\beta_{n-1}]} \]
and from~\cite[1.6]{H} we get 
\[{\HH}^i(A(\Sigma_n)) =
\begin{cases}
k \quad &\mbox{if $i=0$}; \\
k^{n-1} \quad &\mbox{if $i=1$}; \\
0\quad &\mbox{otherwise}.
\end{cases}\]
Hence $\dim_k {\HH}^1(A) \geq n-1$ by Corollary \ref{46}.
\end{example}


\begin{thebibliography}{99}

\bibitem{AP} 
{\bibname I. Assem \and J. A. de la Pe\~na}, `The fundamental groups of a triangular algebra', {\em Comm. Algebra} 24 (1996), no. 1, 187--208.

\bibitem{blm} 
{\bibname M. J. Bardzell, A. C. Locateli \and E. N. Marcos},  `On the Hochschild cohomology of truncated cycle algebras', {\em Comm. Algebra} 28 (2000), no. 3, 1615--1639.

\bibitem{B} {\bibname J. C. Bustamante}, `On the fundamental group of a Schurian algebra', {\em Comm. Algebra}  30  (2002),  no. 11, 5307--5329.
 
\bibitem{CE} {\bibname H. Cartan \and S.Eilenberg}, `Homological algebra', {\em Princeton University Press, Princeton, N. J.}, 1956. xv+390 pp. 

\bibitem{c1} {\bibname C. Cibils}, `On the Hochschild cohomology of finite-dimensional algebras', {\em Comm. Algebra} 16 (1988), no. 3, 645--649.

\bibitem{C} {\bibname C. Cibils}, `Cohomology of incidence algebras and simplicial complexes', {\em J. Pure Appl. Algebra} 56 (1989), no. 3, 221--232.

\bibitem{C2}{\bibname C. Cibils}, `Hochschild cohomology algebra of radical square zero algebras', {\em  Algebras and modules, II (Geiranger, 1996),  93--101, CMS Conf. Proc., 24, Amer. Math. Soc., Providence, RI}, 1998.

\bibitem{crs} {\bibname C. Cibils, M. J. Redondo \and M. Saor\' \i n}, `The first cohomology group of the trivial extension of a monomial algebra', {\em  J. Algebra Appl.} 3  (2004),  no. 2, 143--159. 

\bibitem{cs} {\bibname C. Cibils \and M. Saor\'\i n}, `The first cohomology group of  an algebra with coefficients
in a bimodule', {\em J. Algebra} 237 (2001), no. 1, 121--141.

\bibitem{D} {\bibname P. Dr\"axler}, `Completely separating algebras', {\em J. Algebra} 165  (1994),  no. 3, 550--565.

\bibitem{GR} {\bibname M. A. Gatica \and M. J. Redondo}, `Hochschild cohomology and fundamental groups of incidence algebras', {\em Comm. Algebra}  29  (2001),  no. 5, 2269--2283. 

\bibitem{GR1} {\bibname M. A. Gatica \and M. J. Redondo}, `Hochschild cohomology of incidence algebras as one-point extensions', {\em Linear Algebra Appl.}  365  (2003), 169--181. 
 
\bibitem{GS} {\bibname M. Gerstenhaber \and S. D. Schack}, `Simplicial cohomology is Hochschild cohomology', {\em  J. Pure Appl. Algebra}  30  (1983),  no. 2, 143--156.

\bibitem{H} {\bibname D. Happel}, `Hochschild cohomology of finite--dimen\-sio\-nal algebras', S\'eminaire d'alg\`ebre Paul Dubreuil et
Marie--Paule Malliavin, {\em Lect. Notes Math.} 1404, 108--126, 1989.

\bibitem{IZ} {\bibname K. Igusa \and D. Zacharia}, `On the cohomology of incidence algebras of partially ordered sets', {\em Comm. Algebra} 18 (1990), no. 3, 873-887.

\bibitem{PS} {\bibname J. A. de la Pe\~na \and M. Saor\'\i n}, `On the first Hochschild cohomology group of an algebra', {\em Manuscripta Math.} 104 (2001), 431--442.

\bibitem{R} {\bibname M. J. Redondo}, `Hochschild cohomology: some methods for computations', {\em IX Algebra Meeting USP/UNICAMP/UNESP (S\~ao Pedro, 2001).  Resenhas} 5  (2001),  no. 2, 113--137.

\bibitem{Re1} {\bibname E. Reynaud}, `Algebraic fundamental group and simplicial complexes', {\em J. Pure Appl. Algebra}  177  (2003),  no. 2, 203--214.

\bibitem{Re2} {\bibname E. Reynaud}, `Incidence algebras and algebraic fundamental group', Th\'eories d'homologie, repr\'esentations et alg\`ebres de Hopf. {\em AMA Algebra Montp. Announc.} 2003, Paper 7, 6 pp. (electronic).

\end{thebibliography}
\end{document}